\documentclass[]{amsart}
\usepackage{graphicx}
\usepackage[T1]{fontenc}
\usepackage[latin1]{inputenc}
\usepackage{enumerate}
\usepackage[centertags]{amsmath}
\usepackage{pstricks}
\usepackage{latexsym,amssymb}
\newcommand\car{{\mathbf 1}}
\newcommand\cro[1]{\langle #1 \rangle}
\newcommand\N{{\mathbb N}}

\newcommand\M{{\mathcal M}}
\newcommand\C{{\mathcal C}}

\newcommand{\A}{{\mathfrak  A}}

\newcommand{\R}{{\mathbb R}}
\newcommand{\F}{{\mathcal F}}

\newcommand\Mf{{\mathbf M_f^+}}

\newcommand\espp[1]{{\mathbf E}^0\left[#1\right]}
\newcommand\prp[1]{{\mathbf P}^0\left[#1\right]}

\newcommand\pp{{\mathbf P^0}}
\def\N{{\mathbb N}}
\def\Z{{\mathbb Z}}

\newcommand\Rp{{\R^*_+}}

\newcommand\tend{{\underset{n \rightarrow \infty}{\longrightarrow}}}

\newcommand\proc[1]{\left(#1_t\right)_{t \ge 0}}
\newcommand\procz[1]{\left(#1_t\right)_{t \in \R}}

\newcommand\suite[1]{\left\{#1_n\right\}_{n \in \N}}
\newcommand\suiten[1]{\left\{#1\right\}_{n \in \N}}
\newcommand\suitez[1]{\left\{#1_n\right\}_{n \in \Z}}

{\bf}{\it}
\newtheorem{theorem}{Theorem}

\newtheorem{lemma}{Lemma}

{\bf}{\it}
{\bf}{\it}
\newtheorem{corollary}{Corollary}{\bf}{\it}
\begin{document}

\title{Stability of a processor sharing queue with varying throughput%\thanks{Grants or other notes
}

%\maketitle
\author{P. Moyal}
\address{Laboratoire de Math\'ematiques Appliqu\'ees de Compi\`egne\\
Universit\'e de Technologie de Compiègne\\
D\'epartement G\'enie Informatique\\
Centre de Recherches de Royallieu\\
BP 20 529\\
60 205 COMPIEGNE Cedex\\
FRANCE}
\email{Pascal.Moyal@utc.fr} % Your postal address goes here.

\begin{abstract}
In this paper, we present a stability criterion for Processor Sharing queues, in which the throughput may depend on the number 
of customers in the system (in such cases such as interferences between the users). Such a system is represented by a point measure-valued 
stochastic recursion keeping track of the remaining processing times of the customers.
\end{abstract}

 % insert keywords separated by a semicolon

\maketitle
\section{Introduction}
\label{sec:INTRO}
In this paper, we address the question of stationarity in the general ergodic framework for processor sharing queues, in which the throughput (i.e. the quantity of work achieved by the server(s) per unit of time) may depend on the state of system. More precisely, we assume hereafter that the server(s) 
(it will be clear in the sequel that the effective number of servers does not really matter, only does the quantity of work consumed per unit of time) 
process(es) all the jobs present in the system simultaneously and fairly. Whenever there are $n$ customers in the system, each of them is thus served at 
a rate that depend on $n$, say $r(n)$. The classical case is when $r(n)=1/n$, $n\ge 1$, so that the total throughput equals 
$n.r(n)=1$ whenever the system is non-empty: this is the classical Processor Sharing queue. Hereafter we consider a more general context, 
in which the total throughput may decrease with the number of customers in the system (hence $n.r(n)\le 1$). This is the case for instance in a  wireless network in which the number of users being currently active may decrease the efficiency of the resources. Another case, is when 
the value of $n$ the number of customers does \emph{not} change the nominal service rate $r(n)$, say $r(n)=1$ for all $n$. This corresponds to the classical queue with infinitely many servers.

In both cases and under general stationary ergodic assumptions, Loynes' stability result does not hold since this is not a proper G/G/1 queue 
(the throughput may be less, or larger than one). We address the question of the existence of a stationary version of such queues by representing them with point measure-valued stochastic recursions in the Palm setting, so as to take into account the dependency in the number of customers. This point measures keep track of all the remaining service times of all the customers in the system. Then, it is possible to provide conditions for the existence 
of a stationary version of this sequence, that allow to explicitly construct stationary queues under these assumptions.

This paper is organized as follows. After some preliminaries in section \ref{sec:ORDER}, we present the queueing models we consider in section \ref{sec:PS}. In section \ref{sec:PUR} we study the particular case of the G/G/$\infty$ queue, and in section \ref{sec:PSr} we present a stability 
criterion for generalized processor queues with state-dependent throughput. 
\section{Preliminaries}
\label{sec:ORDER} 
Let $\mathbf{M}_f^+$ and $\C_b$ denote
respectively the set of positive finite measures on $\Rp$ and the set
of bounded continuous functions from $\R$ to $\R$. Equipped with the
\emph{weak topology} $\sigma\left(\mathbf M_f^+,\C_b\right)$,
$\mathbf M_f^+$ is Polish (see \cite{Bill68}). Let $\mathbf 0$ be the zero measure on $\R$ (i.e., such that $\mathbf 0(\mathfrak B)=0$
for any Borel set $\mathfrak B$ on $\R$). 
For any $\mu \in \Mf$ %with compact support, we denote $m(\mu):=\sup\left\{a \in \Rp, \mu\left([a,\infty)\right)>0\right\}$ 
and any measurable $f~:\R \rightarrow \R$, we classically write 
$\cro{\mu,f}:=\int f \,d\mu$. 
%We denote $\C\left(\Mf\right)$, the space of continuous mappings from $\Mf$ into itself. 
Let us denote for any $y \in \R$ and any measurable $f: \R \rightarrow \R$, 
$\tau_yf(.)=f(.-y)\car_{\left\{.>y\right\}}$. Then, for any $\mu \in \Mf $, $\tau_y\mu$ denotes the only element of $\Mf$ s.t. 
$\cro{\tau_y\mu,f}=\cro{\mu,\tau_yf}$.

Let the set $\mathbf M_f^+$ be endowed with the \emph{increasing partial integral order} $\preceq $~: for any two 
$\mu, \nu \in \mathbf M_f^+$, $\mu \preceq \nu$ if $\cro{\mu,f} \le \cro{\nu,f}$ for any measurable non-decreasing function $f$ such that these integrals exist.  
Of course, $\mathbf 0 \preceq \mu$ for any $\mu \in \Mf$. Furthermore, let us remark that
\begin{lemma}\label{lemma:base}
Any sequence of $\Mf$ that is $\preceq$-increasing and bounded above converges for the weak topology. 
\end{lemma}

\begin{proof}
Let $\suite{\mu}$ be a $\preceq $-increasing sequence of $\Mf$ that is bounded above by $\mu \in \Mf$. Then, as easily seen the sequence of non-increasing real functions $\left\{\mu_n\left([.,\infty)\right)\right\}_{n\in\N}$ tends pointwise, and hence (this is Diniz Theorem), uniformly to a non-increasing real function $f$ that is right continuous and has a countable number of discontinuities. Moreover $f(0)\le \mu(\Rp)<\infty$, and we can fully characterize a measure 
$\mu^* \in \Mf$ setting $\mu^*\left((0,x)\right)=f(0)-f(x)$ for all $x\in\Rp$. 
In particular, $\sup_{x\in\Rp}\left|\mu^n\left((0,x)\right)-\mu^*\left((0,x)\right)\right|\tend 0,$ hence $\mu^n$ tends to $\mu^*$ in total variation, which completes the proof.  
\end{proof}

Let now $\M \subset \Mf$ be the subset of finite (simple) counting measures on $\Rp$. 
Any $\mu \in \M \backslash \{\mathbf 0\}$ reads 
$\mu=\sum_{i=1}^{N(\mu)}\delta_{\alpha_i(\mu)},$
where $N(\mu):=\mu(\Rp)$ is the number of atoms of $\mu$, $\delta_{x}$ is the Dirac measure at $x \in \R+$ and
$\alpha_1(\mu)<\alpha_2(\mu)<...<\alpha_{N(\mu)}(\mu)$. Then,  
$\tau_y(\mu)=\sum_{i=1}^{N(\mu)}\delta_{\alpha_i(\mu)-y}\car_{\left\{\alpha_i(\mu)>y\right\}}$
%and for any
%$f\in \C_b$,
%$$\cro{\mu,f}=\sum_{i=1}^{N(\mu)}f\left(\alpha_i(\mu)\right).$$
and for any two $\mu,\nu \in \M  \backslash \left\{\mathbf 0\right\}$, $\mu \preceq \nu$ whenever
\[\left\{\begin{array}{ll} 
(i)\displaystyle &N(\mu) \le N(\nu),\\\\ (ii)\displaystyle& \mbox{for all }i=0,...,N(\mu)-1, \alpha_{N(\mu)-i}(\mu) \le \alpha_{N(\nu)-i}(\nu).
\end{array}\right.\]
We denote for any $\mu \in \M \backslash \left\{\mathbf 0\right\}$, $Z(\mu)=\alpha_{N(\mu)}(\mu)$, the largest atom of $\mu$. Finally, we write $x^+=\max(x,0)$ for any real number $x$, $\sum_{i=j}^k.\equiv 0$ whenever $k<j$ and $\max\left\{\emptyset\right\}\equiv 0$. 

\section{The model}
\label{sec:PS}
Let us first introduce our definitions and assumptions on the queueing systems we shall consider in the sequel. 
Let $\left(\Omega,\F,\mathbf P,\theta_t\right)$ be a probability space furnished with a bijective flow $\proc{\theta}$, under which $\mathbf P$ is stationary and ergodic. Define on $\Omega$ the $\theta_t$-compatible simple point process $\procz{A}$ of points 
$...<T_{-2}<T_{-1}<T_0\le 0< T_1<T_2<...$, that represent the arrival times of the customers in a queue without  
buffer. The process $\procz{A}$ is marked by a sequence $\suitez{\sigma}$, where for all $n \in \Z$, $\sigma_n$ is the service duration requested 
by the customer $C_n$ arrived at time $T_n$.  Also denote for all $n \in \Z$, $\xi_n=T_{n+1}-T_n$, and suppose that the generic r.v. $\sigma$ and 
$\xi$ are integrable. 
We consider that the server(s) follow a generalized Processor Sharing discipline.  By that, we mean that all present customers are taken care of simultaneously, at a rate $r$ that is equal for all customers. An example is of course provided by the classical Processor Sharing queue, but it will be shown in the subsequent sections that significant results can be obtained as well for a wider class of systems. Indeed, it is plausible to assume in many cases, that the amount of work in the 
system might affect the throughput, considering for instance the working cost induced by the switching mechanism in the processor, or the  interferences between the users of a wireless network. In both cases, it is then natural to assume that the rate $r$ is a non-increasing function of the service profile, \emph{i.e.} $\mu \preceq \nu$ implies $r(\mu)\ge r(\nu)$. Hereafter, for the sake of simplicity, we will restrict to the sub-case, where $r$ is a non-increasing function of the number of customers in the system, although it should be clear that all the results below hold as well when $r$ is function of the whole service profile. In other words, at any  $t$, denoting 
$Q(t)$ the number of customers in the system at $t$, each customer is allocated a quantity of work $r(Q_t)$ per unit of time, that is such that 
$r(i) \ge r(j)$ for all $i,j \in \N^*$ such that $i\le j$. Let us illustrate through a naive example the effect of a large number 
of customers on the throughput. 
\begin{center}
\begin{tabular}{ccc}
& \textbf{Nominal service rate} & \textbf{Troughput} \\
\textbf{1 customer} & 1 & 1  \\
\textbf{2 customer}& 0.495 & 0.99 \\
\textbf{3 customer}& 0.3 & 0.9 \\
...& ...& ... \\
\textbf{100 customers}& 0.008 & 0.8  \\
\end{tabular}
%\label{exempleprob}
\end{center}

Provided that $C_n$ is in the system at  $t$, his \emph{remaining processing time} at this instant is the time before his service completion. 
The \emph{service profile} of the system at $t$ is the $\M$-valued process keeping track of the remaining processing times of all the customers 
in the system at $t$:
$$\mu(t)=\sum_{i=1}^{Q(t)}\delta_{\alpha_i(\mu(t))}$$
where $\alpha_1(\mu(t)) \le \alpha_2(\mu(t)) \le .....\le \alpha_{Q(t)}(\mu(t))$ denote the remaining processing times 
of the customers in the system at $t$, ranked in decreasing order. Let $W(t)$ denote the workload at $t$. 
Then, the workload and the congestion processes can easily be recovered from the service profile process by writing for all $t$
\[\begin{array}{ll}
\displaystyle Q(t)&=N(\mu(t)),\\
\displaystyle W(t)&=\cro{\mu(t),I},
\end{array}\]
where $I$ is the identity function. 
The processes $\mu$, $Q$ and $W$ have RCLL paths, and we denote for all $t$ $\mu(t-)=\lim_{s\uparrow \uparrow t} \mu(s)$ (and accordingly, $Q(t-)$ and $W(t-)$). We denote for all $n \in \N$, $\mu_n=\mu(T_n-)$ (respectively $Q_n=Q(T_n-)$, $W_n=X(T_n-)$) the service profile (resp. congestion, workload) just before the arrival of customer $C_n$.

Let $\left(\Omega,\F,\pp\right)$ be the Palm space of $A$, denote
$\theta:=\theta_{T_1}$, $\theta^{-1}$ his measurable inverse and for all $n \in \Z$,
$\theta^n=\theta\circ\theta\circ...\circ\theta\,\,\mbox{ and }\,\,\theta^{-n}=\theta^{-1}\circ\theta^{-1}\circ...\circ\theta^{-1}.$  
Note, that $\pp$ is stationary and ergodic under $\theta$, 
\emph{i.e.}  
for all $\A \in \F$, $\prp{\theta^{-1}\A}=\prp{\A}$ and $\theta\A=\A$ implies 
$\prp{\mathfrak A}=0$ or $1$, and that all $\theta$-contracting event (such that $\prp{\A^c\cap\theta^{-1}\mathfrak A}=0$) is $\theta$-invariant. Denoting $\xi:=\xi_0$ and $\sigma:=\sigma_0$, we have
for all $n\in\Z$,  
$\xi_n:=\xi\circ\theta^n$ and $\sigma_n:=\sigma\circ\theta^n.$

We say that the $E$-valued random sequence $\suite{X}$ is a stochastically recursive sequence (SRS) 
whenever for some random mapping $\phi:\,E\rightarrow E$, 
$$X_{n+1}=\phi\circ\theta^n\left(X_n\right),\, n\in\N,\,\pp-\mbox{ a.s.}.$$ 
For any $E$-valued r.v. $Y$, we denote $\suite{X^{[Y]}}$ the SRS $\suite{X}$ such that $X^{[Y]}_0=Y$, $\pp$-a.s.. 
We follow the formalism of \cite{BacBre02} and formulate the question of stationarity for the SRS $\suite{X}$ in the following terms. 
There exists a \emph{stationary version} of $\suite{X}$ whenever for some $Y$ and for all $n$, $X_n^{[Y]}=Y\circ\theta^n$, $\pp$-a.s., or in other words, provided that the equation 
$$Y\circ\theta=\phi(Y)$$
admits a solution that is a $E$-valued r.v..  
We say that two sequences of r.v. $\suite{X}$ and $\suite{Y}$ couple provided that 
$$\prp{\exists N(\omega), X_n(\omega)=Y_n(\omega) \mbox{ for all }n \ge N(\omega)}=1,$$ 
and that there is strong backwards from $\suite{X}$ with the stationary sequence $\left\{Y\circ\theta^n\right\}$ whenever 
$$\prp{\exists N'(\omega), X_n\circ\theta^{-n}(\omega)=Y(\omega) \mbox{ for all }n \ge N'(\omega)}=1.$$

\begin{lemma}
\label{lemma:recurr}
The sequence $\suite{\mu}$ is stochastically recursive for any rate function $r$: letting for all $\mu \in \M$ and $x \in \Rp$,
\begin{itemize}
\item For all $i \le N(\mu)$,
$$\gamma^r_i(\mu,x)=r(N(\mu)-i+1)\left(x-\sum_{j=1}^{i-1} \alpha_j(\mu)\left(\frac{1}{r(N(\mu)-j+1)}-\frac{1}{r(N(\mu)-j)}\right)\right),$$
\item $i^r(\mu,x)=\max\Biggl\{i \le N(\mu); \alpha_i(\mu)\le\gamma^r_i(\mu,x)\Biggl\},$ % \mbox{ (or 0 if the set is empty)},$$
\item $\gamma^r(\mu,x):=\gamma^r_{(i^r(\mu,x)+1)\wedge 1}(\mu,x),$
\item $\Phi^{r}(\mu,x)=\tau_{\gamma^r(\mu,x)}\mu,$ 
\end{itemize}
we have for any initial profile $\mu_0$ anf for all $n \in \N$, 
% $$\varsigma _{n+1}= \Psi _{\xi_n}\left(\varsigma_n+\delta_{\sigma_n}\right),$$ 
% where for all $x \ge 0$ and all $\mu \in \Mf$, 
% $$\Psi_x(\mu)=\sum_{i=i^r+1}^{N(\mu)}\delta_{\alpha_i(\mu)-\gamma_{i^r}(\mu,x)},$$
\begin{equation}
\label{eq:recurr}
\mu_{n+1}=\Phi^{r}\left(\mu_n+\delta_{\sigma_n},\xi_n\right).
\end{equation}
\end{lemma}

\begin{proof}
%The notations are that of the proof of Lemma \ref{lemma:recurSRPT}. 
Just after the arrival of $C_n$, the service profile reads $\mu:=\mu_n+\delta_{\sigma_n}$. Set $T^{\prime}_{0}:=T_n$ and $\alpha_0\left(\mu\right)=0$. For any $i \in \left\{1,...,N\left(\mu\right)\right\}$, let $T^{\prime}_i$ be the theoretical departure of the customer $\tilde C_i$ whose remaining service  
time at $T_n$ is $\alpha_i(\mu)$. 
%$\tilde C_i$ and $\tilde C_{i-1}$ receive the exact same quantity of service in the interval of time $[T_n,T^{\prime}_i]$, i.e. $\mu^i$. Thus 
The remaining service time of $\tilde C_{i}$ at $T^{\prime}_{i-1}$ is $\alpha_i(\mu)-\alpha_{i-1}(\mu)$, and between $T^{\prime}_{i-1}$ and $T^{\prime}_{i}$, $\tilde C_{i}$ is served at rate $r(N(\mu)-i+1)$. 
%Hence it is easily seen by backwards induction that for all $i \in \left\{1,...,N(\mu)\right\}$, 
Hence we have the induction formula
%%$$T^{\prime}_{i-1}=T^{\prime}_i+\frac{\mu^{i-1}-\mu^i}{r(i-1)},\,i \in \left\{1,...,N(\mu)\right\}.$$
\begin{equation}
\label{eq:inducformula}
T^{\prime}_{i}=T^{\prime}_{i-1}+\frac{\alpha_i(\mu)-\alpha_{i-1}(\mu)}{r(N(\mu)-i+1)},\,i \in \left\{1,...,N(\mu)\right\},
\end{equation}
from which we deduce that for all $i \in \left\{1,...,N(\mu)\right\}$, 
% From this, we deduce for all $i \in \left\{1,...,N(\mu)\right\}$ that 
\begin{equation}
\label{eq:developpformula}
T^{\prime}_i=T_n+\frac{\alpha_i(\mu)}{r(N(\mu)-i+1)}+\sum_{j=1}^{i-1} \alpha_j(\mu)\left(\frac{1}{r(N(\mu)-j+1)}-\frac{1}{r(N(\mu)-j)}\right).
\end{equation} 
For any $i$, customer $\tilde C_i$ leaves the system before $T_{n+1}$ provided that $T^{\prime}_i-T_n\le \xi_n$, which is equivalent to 
$\alpha_i(\mu)\le \gamma^r_i(\mu,\xi_n)$ in view of (\ref{eq:developpformula}).  
% $$\alpha_i(\mu)\le r(N(\mu)-i+1)\left(\xi_n-\sum_{j=1}^{i-1} \alpha_j(\mu)\left(\frac{1}{r(N(\mu)-j+1)}-\frac{1}{r(N(\mu)-j)}\right)\right)=\gamma^r_i(\mu,\xi_n).$$
%%=\gamma_i\left(\varsigma_n+\delta_{\sigma_n},\xi_n\right).$$ 
In particular, $i^r(\mu,\xi_n)$ denotes the index of the last customer leaving the system before $T_{n+1}$ (or $0$ if there is no departure between $T_n$ and $T_{n+1}$). Then the system is not empty at $T_{n+1}-$ provided that $i^r(\mu,\xi_n)<N(\mu)$, and in this case, $\left\{\tilde C_i, i \in \left\{i^r(\mu,\xi_n)+1,N(\mu)\right\}\right\}$ is the set of customers present in the system at $T_{n+1}-$. For such $i>i^r(\mu,\xi_n)$, the remaining service time of $\tilde C_i$ at $T_{n+1}$ is given by 
% his remaining service time at $T^{\prime}_{i^r}$ minus the quantity of work received since $T^{\prime}_{i^r}$, i.e. 
\begin{equation*}
%\label{eq:PS2ter}
\alpha_i(\mu)-\alpha_{i^r(\mu,\xi_n)}(\mu)-r(N(\mu)-i^r(\mu,\xi_n))\left(T_{n+1}-T^{\prime}_{i^r(\mu,\xi_n)}\right)= \alpha_i(\mu)-\gamma^r\left(\mu,\xi_n\right).
\end{equation*} 
Thus the functional mapping the profile at $T_n$ onto the profile at $T_{n+1}-$ reads
$$\Phi^{r}(.,\xi_n): \mu \longmapsto \sum_{i=i^r(\mu,\xi_n)+1}^{N(\mu)}\delta_{\alpha_i(\mu)-\gamma^r(\mu,\xi_n)}.$$
To obtain the announced result, remark that for any $\mu \in \M$ and $x \in \Rp$, 
for any $i<N(\mu)$ we have that 
$$\gamma_{i+1}^r(\mu,x)-\gamma_{i}^r(\mu,x)=\frac{r\left(N(\mu)-i\right)-r\left(N(\mu)-i+1\right)}{r\left(N(\mu)-i+1\right)}
\left(\gamma^r_i(\mu,x)-\alpha_i(\mu)\right),$$ 
which is nonnegative if and only if $i\le i^r(\mu,x).$ Hence, 
\begin{equation}
\label{eq:PS2}
\gamma^r(\mu,x)=\max_{1\le i \le N(\mu)}\gamma^r_i(\mu,x),
\end{equation}
%
% for all $i\le i^r(\mu,x)$, $\alpha_i(\mu)<\gamma^r_i(\mu,x)$. This implies that $\gamma^r_{i+1}(\mu,x)\ge \gamma^r_i(\mu,x)$ whenever $i<N(\mu)$ as well, as easily derived. 
%% that 
%% \begin{multline*}
%% \gamma_{i-1}(\mu,x)-\gamma_i(\mu,x)=\mu^{i-1}-\mu^i+r(i-1)\left(x-\frac{\mu^{i-1}}{r(i-1)}-\sum_{j=i}^{+\infty}\mu^j\left(\frac{1}{r(j)}-\frac{1}{r(j-1)}\right)\right)\\\shoveright{-r(i)\left(x-\frac{\mu^{i}}{r(i)}-\sum_{j=i+1}^{+\infty}\mu^j\left(\frac{1}{r(j)}-\frac{1}{r(j-1)}\right)\right)}\\
%% \ge \mu^{i-1}-\mu^i+r(i-1)\left(\frac{\mu^{i}}{r(i)}-\frac{\mu^{i-1}}{r(i-1)}-\mu^i\left(\frac{1}{r(i)}-\frac{1}{r(i-1)}\right)\right)=0.
%% \end{multline*}
% Hence we have 
% \begin{equation}
% \label{eq:PS2}
% \gamma^r_j(\mu,x)\le \gamma^r_k(\mu,x) \mbox{ for all $j \le k \le \left(i^r(\mu,x)+1\right)\wedge N(\mu)$},
% \end{equation}
and in particular   
$\Phi^{r}(\mu,\xi_n)=\tau_{\gamma^r(\mu,\xi_n)}\mu,$ $\pp$-a.s..  
\end{proof}
For a fixed $x\in\R+$, the two following monotonicity properties of the mappings $\Phi^r(.,x)$ hold, as shown in Appendix. 
\begin{lemma}
\label{lemma:monotonephi}
For any $x\in\R+$ and any rate function $r$, the mapping $\Phi^r(.,x)$ is non-decreasing from $\M$ into itself.
\end{lemma}
\begin{lemma}
\label{lemma:monotoner}
For any $x\in\R+$ and any $\mu \in\M$, for any two rate functions $r$ and $\tilde r$ such that 
$r(i)\le \tilde r(i)$ for all $i \in \N^*$, $\Phi^r(\mu,x) \succeq  \Phi^{\tilde r}(\mu,x).$ 
\end{lemma}

\section{The pure delay system}
\label{sec:PUR}
Let us first consider the case, where the rate function is constant with respect to the size of the system, say $r(i)=1$ for any $i \ge 1$. 
This corresponds to the classical "pure delay" G/G/$\infty$ queue: all present customers are simultaneously served at unit rate, and hence spend in the 
system a time equal to their service duration, which is equivalent to say that there is an infinity a servers. 
In this case, the recursive equation (\ref{eq:recurr}) driving the service profile sequence (for which a diffusion approximation is given in \cite{moyal06_2} in the M/GI/$\infty$ case)  specializes to 
\begin{equation}
\label{eq:recurPD}
\mu_{n+1}=\tau_{\xi_n}\left(\mu_n+\delta_{\sigma_n}\right) 
\end{equation}
and a stationary service profile for the queue is a solution to the equation 
\begin{equation}
\label{eq:recurstatPD}
\mu\circ\theta=\tau_{\xi}\left(\mu+\delta_{\sigma}\right). 
%\mu\circ\theta=\left(\mathcal R\left(\mu+\sigma.\mathbf {e_{N(\mu)+1}}\right)-\xi.\mathbf 1\right)^+.
\end{equation}
The following lemma (see \cite{moyal07}) will be used in the sequel.  
\begin{lemma}
\label{lemma:solrecurL} 
There exists a unique $\pp$-a.s. finite solution to the equation 
\begin{equation}
\label{eq:recurstatL}
L\circ\theta=\left[\max\left\{L,\sigma\right\}-\xi\right]^+,
\end{equation} 
given by
\begin{equation}
  \label{eq:defL}
  L:=\left[\sup_{j \in \N^*} \left(\sigma_{-j}-\sum_{i=1}^{j}\xi_{-i}\right)\right]^+.
\end{equation}
\end{lemma}

\begin{proof}
\emph{Existence.}    
 Loynes' Theorem for stochastic recurrences   
(see \cite{Loynes62}, \cite{BacBre02}) can be applied since the mapping 
$x\mapsto\left[\max\left\{x,\sigma\right\}-\xi\right]^+$ 
is $\pp$-a.s. continuous and non-decreasing. The minimal solution 
$L$ to (\ref{eq:recurstatL}) classically reads as the $\pp$-a.s. limit of Loynes's sequence  
$\suiten{L^{[0]}_n\circ\theta^{-n}}$, where $\suite{L^{[0]}}$ is the initially null SRS that is defined by  
$$L^{[0]}_{n+1}=\left[\max\left\{L^{[0]}_n,\sigma_n\right\}-\xi_n\right]^+\mbox{ for all }n \in \N.$$ 
It is routine to check from Birkhoff's ergodic theorem (and the fact that $\sigma$ is not identically zero) 
that $L$ is $\pp$-a.s. finite.\\ 
\emph{Uniqueness.} 
Let $\tilde L$ be a solution to (\ref{eq:recurstatL}). First, remark that if $\tilde L>\sigma$, 
$\pp$-a.s. would imply that on a $\pp$-a.s. event,
\begin{equation*}
%\label{eq:lemergo1}
\tilde L \circ \theta >0  \Leftrightarrow
\tilde L \circ \theta = \tilde L-\xi, 
\end{equation*} 
a contradiction to the ergodic Lemma. 
Hence in view of the minimality of $L$, we have that $$\prp{\tilde L=L}=\prp{\tilde L\circ\theta \le L\circ\theta} \ge \prp{\tilde L \le \sigma}>0,$$ 
which implies that $\left\{\tilde L=L\right\}$ is $\pp$-almost sure since it is $\theta$-contracting.  
\end{proof} 
We can now state the following result.
\begin{theorem}
\label{thm:recur1}
The equation (\ref{eq:recurstatPD}) admits a finite solution, given by 
$$\mu^{\text{\tiny{PD}}}=\sum_{i=1}^{\infty}\delta_{\left(\sigma_{-i}-\sum_{j=1}^{i}\xi_{-j}\right)}\car_{\left\{\sigma_{-i}\ge\sum_{j=1}^{i}\xi_{-j}\right\}}.$$
Moreover, provided that 
\begin{equation}
\label{eq:solrecur}
\prp{L\le 0}>0, 
\end{equation}
this solution is unique and for all $\zeta$ such that $Z(\zeta)\le L$, $\pp$-a.s, the sequence $\suite{\mu^{[\zeta]}}$ converges with strong backwards coupling to $\mu^{\text{\tiny{PD}}}.$ 
\end{theorem}
\begin{proof}
\emph{Existence.} 
It is a straightforward consequence of Birkhoff's ergodic theorem that  
$$\prp{\mu^{\text{\tiny{PD}}} \in \M}=\prp{\mbox{Card}\,\left\{i \in \N^*,\sigma_{-i}-\sum_{j=1}^{i}\xi_{-j}\ge 0\right\}<\infty}>0.$$ 
This $\theta$-contracting event is thus $\pp$-almost sure. On another hand, in view of Lemma \ref{lemma:monotonephi}, the mapping $\mu \mapsto \tau_{\xi}\left(\mu+\delta_{\sigma}\right)$ is $\pp$-a.s. 
non-decreasing from $\M$ into itself. It is furthermore continuous for the weak topology, as easily checked from the fact that
for any $\M$-valued sequence $\suite{\nu}$ tending weakly to $\nu$, for any $x,s \in \R+$ and any $\phi \in \C_b$, 
$$\cro{\tau_x\nu_n+\delta_{s},\phi}=\int\phi(y-x)\,d\nu_n(y)+\phi(s)\tend \int\phi(y-x)\,d\nu(y)=\cro{\tau_x\nu+\delta_{s},\phi}.$$
 Thus, we can follow the steps of Loynes' construction (Lemma \ref{lemma:base}), to conclude that $\mu^{\text{\tiny{PD}}}$ is the $\preceq$-minimal solution of (\ref{eq:recurstatPD}) since it is the $\pp$-a.s. limit of the sequence given for all $n \in\N$ by   
$$\mu_n^{[\mathbf 0]}\circ\theta^{-n}=\sum_{i=1}^{\infty}\delta_{\left(D_{-i}-\sum_{j=1}^{i}\xi_{-j}\right)}\car_{\left\{D_{-i}\ge\sum_{j=1}^{i}\xi_{-j}\right\}}.$$ 
\emph{Uniqueness.} 
It is easily checked that for any solution $\mu$ of (\ref{eq:recurstatPD}), 
$$Z(\mu)\circ\theta=Z\left(\tau_{\xi}\left(\mu+\delta_{\sigma}\right)\right)=\left[Z(\mu)\vee \sigma-\xi\right]^+,$$ 
hence $Z(\mu)=L$, $\pp$-a.s.. Moreover, since $\mu^{\text{\tiny{PD}}}$ 
is the minimal solution of (\ref{eq:recurstatPD}), we have that 
$$\left\{\mu=\mu^{\text{\tiny{PD}}}\right\} \supseteq
\left\{\mu= \mathbf 0 \right\}=\left\{Z(\mu)=0\right\}=\left\{L=0\right\}.$$ 
Hence, whenever (\ref{eq:solrecur}) holds, the event $\left\{\mu=\mu^{\text{\tiny{PD}}}\right\}$ has a positive probability. Since it is $\theta$-invariant, it is $\pp$-almost sure.

\noindent \emph{Coupling.} %Let us assume that (\ref{eq:solrecur}) holds. 
Let $\zeta$ be a $\M$-valued r.v. such that 
$Z(\zeta)\le L$, $\pp$-a.s.. It is easy to construct another $\M$-valued r.v. $\tilde\zeta$ such that $\zeta \preceq \tilde \zeta$ and $Z(\tilde\zeta)=L$, $\pp$-a.s. by setting e.g. $\tilde \zeta=\sum_{i=1}^{N(\zeta)-1}\delta_i(\zeta)+\delta_{L}$. 
From Lemma \ref{lemma:monotonephi}, it follows by induction 
that $\mu_n^{[\zeta]} \preceq \mu_n^{[\tilde \zeta]}$, $\pp$-a.s. for all $n \in \N$. Remark now that for all $n\in\N$, $Z\left(\mu_n^{[\tilde \zeta]}\right)=L\circ\theta^n$, as easily checked by induction. Hence, for all $n\in\N$, we have 
$$\mathcal E_n:=\left\{L\circ\theta^n=0\right\}=\left\{Z\left(\mu_n^{[\tilde \zeta]}\right)=0\right\}=\left\{\mu_n^{[\tilde \zeta]}=\mathbf 0\right\}\subseteq \left\{\mu_n^{[\zeta]}=\mathbf 0\right\}.$$ 
Therefore, $\suite{\mathcal E}$ is a stationary sequence of renovating events of length 1 for $\suite{\mu^{[\zeta]}}$ (see \cite{MR52:12118,Foss92})  for any $\zeta$ such that $Z(\zeta)\le L$, $\pp$-a.s.. Assumptions (\ref{eq:solrecur}) implies the coupling property for such an initial condition in view of Corollary 2.5.1 of 
\cite{BacBre02}. 
\end{proof}
As simple consequences of the latter result, let us remark the following coupling properties.
\begin{corollary} 
Under condition (\ref{eq:solrecur}), for any 
 $\zeta$ such that $Z(\zeta)\le L$, $\pp$-a.s,
\begin{enumerate}
\item [(i)] $\suite{X^{[N(\zeta)]}}$ converges with strong backwards coupling to $N\left(\mu^{\text{\tiny{PD}}}\right);$ 
\item [(ii)] $\suite{W^{[\cro{\zeta,I}]}}$ converges with strong backwards coupling to $\cro{\mu^{\text{\tiny{PD}}},I}.$
%\item [(iv)] For any $i \in \N^*$, $\alpha_i(\mu_n^{\mathbf 0 }) \convloi \alpha_i(\kappa).$
\end{enumerate}
\end{corollary}

\section{Processor Sharing queues}
\label{sec:PSr}
We shall now consider the case, where the rate function depends on the number of customers in the system at current time. 
We assume hereafter that the non-decreasing function $r$ is such that
\begin{equation}
\label{eq:hypor1}
\sup_{n\in \N^*} n.r(n)\le 1,
\end{equation} 
\begin{equation}
\label{eq:hypor2}
K_r=\inf_{n\in \N^*} n.r(n)>0.%\mbox{ For all }i \in \N^*,  K_r\le ir(i) \le 1\mbox{ for some constant }K_r>0.
\end{equation}
Assumption (\ref{eq:hypor1}) amounts to say that there is a single server, in that the throughput at time $t$, given by $Q(t).r(Q(t))$, may not exceed one. 
A typical case is the classical Processor Sharing queue: assume that $r(n)=n^{-1}$ for any $n$ (and hence $K_r=1$), meaning that all customers 
are served at a rate that is inversely proportional to the number of customers. In that case the server works at unit rate whatever the number of 
customers in the system. Whenever $K_r<1$, the number of customers affects the velocity of service, so that the total throughput may be less than one. We assume nevertheless in (\ref{eq:hypor2}) that a minimal throughput $K_r$ is granted for a given $r$, \emph{i.e.} the server always achieves at least $K_r$ unit of work 
per unit of time. An example is provided by the following 
idealistic scenario: the server works at unit rate whenever there is only one customer in the system ($r(1)=1$), and the interferences (or 
operating cost) when there are several 
customers in service at the same time decreases by half the efficiency of the server, so that $r(i)=1/(2i)$ for any $i \ge 2$, which implies in particular that (\ref{eq:hypor2}) is met for $K_r=1/2$.
  
In view of Lemma \ref{lemma:recurr}, a stationary service profile is a solution to the equation 
\begin{equation}
\label{eq:recurstatr}
\mu\circ\theta=\Phi^{r}\left(\mu +\delta_{\sigma},\xi\right).  
\end{equation} 
%where the random mapping $\Phi^{r}(.,\xi)$ is defined in Lemma \ref{lemma:recurr}. 
We have the following result. 
\begin{theorem}
Let $r$ be a rate function satisfying assumptions (\ref{eq:hypor1}) and (\ref{eq:hypor2}). Then provided that 
\begin{equation}
\label{eq:condr}
\espp{\sigma}<K_r\espp{\xi},
\end{equation}
the equation (\ref{eq:recurstatr}) admits a unique finite solution $\mu^{r}$. Moreover, for any $\M$-valued r.v. $\zeta$ such that 
$\cro{\zeta,I}\le W^{K_r}$, $\pp$-a.s. (where $W^{K_r}$ is the unique solution of (\ref{eq:lindleystatKr})), 
the sequence $\suite{\mu^{[\zeta]}}$ converges with strong backward coupling to $\mu^{r}$. 
\end{theorem}
\begin{proof}
\emph{Existence.} 
Fix $r$ satisfying (\ref{eq:hypor1}) and (\ref{eq:hypor2}). From Loynes's fundamental stability result, the equation 
\begin{equation}
\label{eq:lindleystatKr}
W\circ\theta=\left[W+\sigma-K_r\xi\right]^+
\end{equation}
admits a unique $\pp$-a.s. finite solution, say $W^{K_r}$, provided that (\ref{eq:condr}) holds. 
Let $\tilde r$ be the rate function such that for all $\mu \in \M$, 
$\tilde r(\mu)=K_r/N(\mu)$, so that 
the throughput under $\tilde r$ always equals $K_r$ whenever the system is non-empty. 
Let $\zeta$ be a $\M$-valued r.v. such that $\cro{\zeta,I} \le W^{K_r}$ and $$\tilde \zeta=\zeta+\delta_{W^{K_r}-\cro{\zeta,I}}\car_{W^{K_r}>\cro{\zeta,I}}.$$ Is then clear 
that $\cro{\tilde\zeta,I}=W^{K_r}$. Moreover, we have $\pp$-a.s. for all $n\in\N$ 
$$\cro{\mu^{\tilde r,[\tilde \zeta]}_{n+1},I}=\left[\cro{\mu^{\tilde r,[\tilde \zeta]}_{n},I}+\sigma_n-K_r\xi_n\right]^+,$$ 
as the throughput equals $K_r$ at any time (as easily checked from Lemma \ref{lemma:recurr}), so that 
$\cro{\mu^{\tilde r,[\tilde \zeta]}_{n+1},I}=W^{K_r}\circ\theta^n$ for all $n\in\N$. 
On another hand, $\zeta \preceq \tilde \zeta$, hence in view of Lemmas \ref{lemma:monotonephi} and \ref{lemma:monotoner}, 
an immediate induction shows that  
$\mu^{r,[\zeta]}_{n}\preceq \mu^{\tilde r,[\tilde\zeta]}_{n}$ for all $n\in\N$, which implies in turn that  
$$\cro{\mu^{r,[\zeta]}_{n},I}\le \cro{\mu^{\tilde r,[\tilde\zeta]}_{n},I}=W^{K_r}\circ\theta^n \mbox{ for all }n\in\N.$$
% $$\mu^{r,[\zeta]}_{n+1}=\Phi^r\left(\mu^{r,[\zeta]}_{n}+\delta_{\sigma_n},\xi_n\right) \preceq \mu^{\tilde r,[\tilde \zeta]} $$
% Denote now for all $\zeta \in \M$ and all $n\in\N$, $\mu_n^{r,[\zeta]}$ the service profile at $T_n-$ given 
% the service profile at $T_0-$ equals $\zeta$, $\pp$-a.s..  
%%%, $\mu^{r,[\mathbf 0]}_n$ the service profile at $T_n-$ whenever the system is empty at $T_0-$, and 
%%$$W^{r,[0]}_n=\sum_{i \in \N^*}(\mu^{r,[0]}_n)^i$$
%%the workload in the system at $T_n-$ given the system is empty at $T_0-$. 
% In view of (\ref{eq:hypor2}), service 
% is consumed at a rate larger than $K_r$ whenever the system is non-empty, thus we have 
% $$\cro{\mu^{r,[\zeta]}_{n+1},I} \le \left[\cro{\mu^{r,[\zeta]}_{n},I}+\sigma_n-K_r\xi_n\right]^+,\,n \in \N.$$ 
% Hence an immediate induction shows that whenever $\zeta \le W^{K_r}$, $\pp$-a.s., for all $n\in\N$, 
% $$\cro{\mu^{r,[\zeta]}_{n},I} \le W^{K_r}\circ\theta^n,\,\pp-\mbox{a.s.}.$$ 
Therefore, for all $n \in \N$, on $\mathfrak A_n:=\left\{W^{K_r}\circ\theta^n=0\right\}$, we have that 
$\cro{\mu^{r,[\zeta]}_{n},I}=0$, hence $\mu^{r,[\zeta]}_n=\mathbf 0$ and 
$$\mu^{r,[\zeta]}_{n+1}=\Phi^{r}\left(\delta_{\sigma_n},\xi_n\right).$$  
Therefore $\suite{\mu^{r,[\zeta]}}$ admits $\suite{\mathfrak A}$ as a stationary sequence of renovating events of length 1. Furthermore, 
the event $\mathfrak A_0=\left\{W^{K_r}=0\right\}$ has a striclty positive probability since the contrary would imply that 
$$\espp{W^{K_r}\circ\theta-W^{K_r}}=\espp{\sigma-K_r\xi}<0,$$ 
an absurdity in view of the ergodic Lemma. Then it follows from \cite{BacBre02}, Th. 2.5.3, that there is strong backwards coupling of 
$\mu^{r,[\zeta]}_n$ with the stationary sequence $\left\{\mu^{r}\circ\theta^n\right\}_{n\in\N}$, where 
$\mu^r$ is a proper solution to (\ref{eq:recurstatr}). 

\emph{Uniqueness.} 
Fix $r$ and $\tilde r$ be as above. 
There exists a solution $\mu^{\tilde r}$ to (\ref{eq:recurstatr}). Then, we have $\pp$-a.s. 
%$$\cro{\mu^{\tilde r},I}\circ\theta=\cro{\mu^{\tilde r}\circ\theta,I}=\left[\cro{\mu^{\tilde r},I}+\sigma-K_r\xi\right]^+,$$
$$\cro{\mu^{\tilde r},I}\circ\theta=\cro{\Phi^{\tilde r}\left(\mu^{\tilde r}+\delta_{\sigma},\xi\right),I}=\left[\cro{\mu^{\tilde r},I}+\sigma-K_r\xi\right]^+,$$ 
hence $\cro{\mu^{\tilde r},I}$ equals $W^{K_r}$, $\pp$-a.s.. Moreover, on $\left\{\cro{\mu^r,I}\le W^{K^r}\right\}$, we have in view of 
Lemma \ref{lemma:recurr} that
$$\cro{\mu^r,I}\circ\theta\le \cro{\Phi^{\tilde r}\left(\mu^r+\delta_{\sigma},\xi\right),I}=\left[\cro{\mu^r,I}+\sigma-K_r\xi\right]^+\le W^{K_r}\circ\theta,\,\,\pp-\mbox{a.s.},$$ 
thus the event $\left\{\cro{\mu^{r},I}\le W^{K_r}\right\}$ is $\theta$-contracting. Moreover, 
 $$\prp{\cro{\mu^r,I}\le W^{K_r}}\ge \prp{\cro{\mu^r,I}=0}>0,$$ 
as another consequence of (\ref{eq:condr}) and the ergodic Lemma. 
Therefore, $\cro{\mu^r,I}\le W^{K_r}$, $\pp$-a.s., so that 
$$\A_n \subseteq  \left\{\cro{\mu^r,I}\circ\theta^n=0\right\}=\left\{\mu^{r}\circ\theta^n=\mathbf 0\right\}.$$ 
Consequently, $\suite{\A}$ is a stationary sequence of renovating events of length 1 for $\left\{\mu^r\circ\theta^n\right\}_{n\in\N}$ for any 
solution $\mu^r$ of the equation (\ref{eq:recurstatr}) associated to the rate $r$. Since $\prp{\A_0}$, there exists a unique solution to (\ref{eq:recurstatr}) in view of Remark 2.5.3. in \cite{BacBre02}.     
\end{proof}
We have in particular:
\begin{corollary} 
Under condition (\ref{eq:condr}), for any 
$\zeta$ such that 
$\cro{\zeta,I}\le W^{K_r}$, $\pp$-a.s., 
\begin{enumerate}
\item [(i)] $\suite{X^{[N(\zeta)]}}$ converges with strong backwards coupling to $N\left(\mu^{r}\right);$ 
\item [(ii)] $\suite{W^{[\cro{\zeta,I}]}}$ converges with strong backwards coupling to $\cro{\mu^{r},I}.$
%\item [(iv)] For any $i \in \N^*$, $\alpha_i(\mu_n^{\mathbf 0 }) \convloi \alpha_i(\kappa).$
\end{enumerate}
\end{corollary}

% The Acknowledgements are an un-numbered section

% Acknowledgements text here

\appendix
%and then carry on using the \section and \subsection commands, as above.

\section{Proofs of monotonicity}
For easy checking, we present hereafter the details of the derivations proving Lemmas \ref{lemma:monotonephi} and \ref{lemma:monotoner}.  

\begin{proof}[Proof of Lemma \ref{lemma:monotonephi}]
We fix again $x\in\R+$ and $\mu,\,\nu \in \M$ such that $\mu \preceq \nu$. Whenever  
$i^r(\mu,x)< N(\mu)$ (otherwise $\Phi^r(\mu,x)=\mathbf 0$), we have that 
\begin{multline*}
\sum_{j=1}^{N(\nu)-N(\mu)+i^r(\mu,x)} \alpha_j(\nu)\left(\frac{1}{r(N(\nu)-j+1)}-\frac{1}{r(N(\nu)-j)}\right)\\ \ge \sum_{j=1}^{i^r(\mu,x)} \alpha_j(\mu)\left(\frac{1}{r(N(\mu)-j+1)}-\frac{1}{r(N(\mu)-j)}\right),
\end{multline*}
which implies that 
\begin{multline*}
\alpha_{N(\nu)-N(\mu)+i^r(\mu,x)+1}(\nu)\ge \alpha_{i^r(\mu,x)+1}(\mu) \\
%%> \frac{\xi-\sum_{i=1}^{i_0(\mu,\xi)+1}\alpha_i(\mu)}{N(\mu)-(i_0(\mu,\xi)+1)}\\
\ge r(N(\mu)-i^r(\mu,x))\left(x-\sum_{j=1}^{i^r(\mu,x)} \alpha_j(\mu)\left(\frac{1}{r(N(\mu)-j+1)}-\frac{1}{r(N(\mu)-j)}\right)\right)\\
%\shoveleft{\ge r(N(\mu)-i^r(\mu,x)).}\\\left(x-\sum_{j=1}^{N(\nu)-N(\mu)+i^r(\mu,x)} \alpha_j(\nu)\left(\frac{1}{r(N(\nu)-j+1)}-\frac{1}{r(N(\nu)-j)}\right)\right).\\
%%\ge \frac{\xi-\sum_{i=1}^{i_0(\mu,\xi)+N(\nu)-N(\mu)+1}\alpha_i(\nu)}{N(\nu)-\left(N(\nu)-N(\mu)+i_0(\mu,\xi)+1\right)}.
\ge \gamma_{N(\nu)-N(\mu)+i^r(\mu,x)+1}^r(\nu,x).
\end{multline*}
This means that $i_0(\nu,x)\le N(\nu)-N(\mu)+i_0(\mu,x),$ i.e. $N\left(\Phi^{r}(\mu,x)\right)\le N\left(\Phi^{r}(\nu,x)\right).$ 
Hence in view of (\ref{eq:PS2}), we have 

\begin{multline*}
\gamma(\mu,\xi)=\gamma^r_{i^r(\mu,\xi)+1}\left(\mu,x\right)\ge \gamma^r_{(i^r(\nu,\xi)+N(\mu)-N(\nu))^++1}(\mu,x)\\
\ge r\left(N(\nu)-i^r(\nu,x)\right)\left(x-\sum_{j=1}^{i^r(\nu,x)}\alpha_j(\nu)\left(\frac{1}{r\left(N(\nu)-j+1\right)}-\frac{1}{r\left(N(\nu)-j\right)}\right)\right)\\=\gamma^r(\nu,x),
\end{multline*}
which clearly implies that $\Phi^r(\mu,x) \preceq \Phi^r(\nu,x)$. 
\end{proof} 

\begin{proof}[Proof of Lemma \ref{lemma:monotoner}]
We now fix $\mu \in \M$ and $x\in\R+$. For any two rate functions $r$ and $\tilde r$ such that 
$r(i)\le \tilde r(i)$ for any $i \in\N^*$, the induction formula (\ref{eq:inducformula}) straightforwardly shows 
that $i^r(\mu,x)\ge i^{\tilde r}(\mu,x) $ i.e. $N\left(\Phi^r(\mu,x)\right) \le N\left(\Phi^{\tilde r}(\mu,x)\right)$. Hence, as in the 
previous proof, $\gamma^r(\mu,x)\le  \gamma^{\tilde r}(\mu,x)$. 
\end{proof} 

\section*{Acknowledgements}
The author would like to warmly thank Bryan Renne for useful discussions in Dublin.

\end{document}